\documentstyle[11pt,amsfonts,fullpage]{amsart}

\def\endproof{\qed\smallskip}
\def\blacksquare{\hbox to .60em{\vrule width .60em height .60em}}

\newtheorem{theorem}{Theorem}[section]

\newtheorem{remark}[theorem]{Remark}

\newcommand{\bbgin}{\begin}

\begin{document}

\title[]{L$^{\bf 2}$ Curvature and Volume Renormalization of AHE 
Metrics on 4-Manifolds}

\author[]{Michael T. Anderson}

\thanks{Partially supported by NSF Grant DMS 0072591}
\thanks{2000 {\it Math. Sci. Classification}. Primary: 53C25, 53C80, 
Secondary: 58J60}

\maketitle

\abstract
This paper relates the boundary term in the Chern-Gauss-Bonnet formula 
on 4-manifolds $M$ with the renormalized volume $V$, as defined in the 
AdS/CFT correspondence, for asymptotically hyperbolic Einstein metrics 
on $M$. In addition we compute and discuss the differential or variation 
$dV$ of $V$, or equivalently the variation of the $L^{2}$ norm of the Weyl 
curvature, on the space of such Einstein metrics.
\endabstract

\setcounter{section}{-1}

\section{Introduction.}

\setcounter{equation}{0}

 The Chern-Gauss-Bonnet formula for a compact Riemannian 4-manifold 
$(M, g)$ without boundary states that
\begin{equation} \label{e0.1}
\frac{1}{8\pi^{2}}\int_{M}(|R|^{2} -  4|z|^{2})dV = 
\frac{1}{8\pi^{2}}\int_{M}(|W|^{2} -  \frac{1}{2}|z|^{2} + 
\frac{1}{24}s^{2})dV = \chi (M), 
\end{equation}
where $R, W, z, s$ are respectively the Riemann, Weyl, trace-free Ricci 
and scalar curvatures.

 In particular, if $g$ is an Einstein metric, then $z =$ 0, and so 
Einstein metrics minimize the $L^{2}$ norm of the curvature over all 
metrics on $M$. Hence the $L^{2}$ norm of the full curvature of an 
Einstein metric on $M$ is apriori bounded by the topology of $M$.

 If $(M, g)$ is a compact 4-manifold with non-empty boundary, then 
(0.1) no longer holds; there is a correction or defect term given by 
certain curvature integrals over the boundary $\partial M.$ If $(M, g)$ 
is complete and open, then the boundary integrals relate to the 
asymptotic geometry of $(M, g)$.

 When $(M, g)$ is a complete non-compact, Ricci-flat 4-manifold, then 
the defect term in (0.1) is easily identified if the manifold 
asymptotically approaches that of a quotient of ${\Bbb R}^{4},$ i.e. 
$M$ is asymptotically locally Euclidean (ALE), flat, (AF) or locally 
flat, (ALF), c.f. [5].

 In this paper, we consider this issue when $(M, g)$ is an Einstein 
4-manifold of negative scalar curvature, which is asymptotically 
hyperbolic. To define this, let $M$ be an arbitrary compact, connected 
and oriented 4-manifold with non-empty boundary $\partial M;$ we do not 
assume that $\partial M$ is connected. According to Penrose, c.f. [14] 
and also [12], a complete metric $g$ on $M$ is {\it  conformally 
compact}  if there is a smooth defining function $\rho $ on $\bar M = 
M\cup\partial M,$ i.e. $\rho (\partial M) =$ 0, $d\rho  \neq $ 0 on 
$\partial M$ and $\rho  > $ 0 on $M$, such that the metric
\begin{equation} \label{e0.2}
\bar g = \rho^{2}\cdot  g, 
\end{equation}
extends to a smooth metric on $\bar M.$ We require that $\bar g$ is at 
least $C^{3}$ smooth up to $\partial M,$ although this condition could 
be relaxed somewhat.

 Conversely, if $\bar g$ is any smooth Riemannian metric on $\bar M$ 
and $\rho $ is any $C^{1}$ defining function, then $g \equiv  \rho^{-2}\cdot 
\bar g$ gives a complete conformally compact metric on the open 
manifold $M$.

 The defining function $\rho $ is not unique, since it can be 
multiplied by any smooth positive function on $\bar M.$ Hence, both the 
metric $\bar g$ and its induced metric $\gamma $ on $\partial M$ are 
not uniquely defined by $(M, g)$. However, the conformal class $[\gamma 
]$ of the metric $\gamma  = \bar g_{T\partial M}$ is uniquely 
determined by the complete metric $g$; $(\partial M, [\gamma ])$ is 
called the {\it conformal infinity} of $(M, g)$. Conversely, any 
conformal class $[\gamma ]$ on $\partial M$ is the conformal infinity 
of a complete metric on $M$.

 When $(M, g)$ is a complete conformally compact Einstein metric with 
$Ric_{g} = - 3g,$ then the sectional curvatures of $g$ necessarily 
approach $- 1$ uniformly at infinity at an exponential rate, c.f. (1.3) 
below or [8]. Such manifolds are called asymptotically hyperbolic (AH).

 The study of complete AH Einstein (AHE) manifolds has become very 
active recently due to the AdS/CFT correspondence in string theory, 
c.f. [16] and references therein. In order to produce an effective 
gravitational action, one needs to renormalize the volume of such a 
metric, since the volume itself is obviously infinite. If $(M, g)$ is 
an AH Einstein 4-manifold and $\rho $ is any defining function on $M$, 
then one has the following asymptotic expansion for the volume of 
compact domains $B(r) = \{log \rho^{-1} \leq  r\}$ in $M$ as $r 
\rightarrow  \infty ;$
\begin{equation} \label{e0.3}
vol B(r) = v_{o}e^{3r} + v_{1}e^{r} + V + o(1). 
\end{equation}
The coefficients $v_{o}$ and $v_{1}$ depend on the geometry of 
$(\partial M, \gamma )$ as well as the defining function $\rho$ 
in this generality.

 Clearly, since there are numerous defining functions, the exact 
exponential growth rates in $r$, as well as the coefficients, depend on 
the compactification $(\bar M, \bar g),$ and are not defined 
intrinsically w.r.t. $(M, g)$. However, the constant term $V$ in (0.3) 
is an invariant of $(M, g)$, i.e. is independent of the choice of $\rho 
.$ This is one of the elementary consequences of the AdS/CFT 
correspondence, c.f. [16]; a proof appears in [7].

 The first purpose of this paper is to relate the renormalized volume 
$V$ in (0.3) with the Chern-Gauss-Bonnet theorem in dimension 4.
\bbgin{theorem} \label{t 0.1.}
  Let $(M, g)$ be a complete AH Einstein 4-manifold. Then, up to a 
constant, the boundary term at infinity in the Chern-Gauss-Bonnet 
formula renormalizes the volume in the sense of (0.3). In fact,
\begin{equation} \label{e0.4}
\frac{1}{8\pi^{2}}\int_{M}|W|^{2}dV = \chi (M) -  \frac{3}{4\pi^{2}}V. 
\end{equation}
\end{theorem}

 An analogous result holds for arbitrary AH metrics on $M$ which are 
suitably asymptotic to an Einstein metric at infinity, c.f. Remark 1.2. 
Of course (0.4) shows directly that $V$ is an intrinsic invariant of 
$(M, g)$, independent of any compactification $(\bar M, \bar g).$

 One thus has the following universal upper bound on $V$ for any AH 
Einstein metric:
\begin{equation} \label{e0.5}
V \leq  \frac{4\pi^{2}}{3}\chi (M), 
\end{equation}
with equality if and only if $(M, g)$ is hyperbolic. Even when $(M, g)$ 
is hyperbolic, i.e. $M = H^{4}(- 1)/\Gamma ,$ (0.5) gives non-trivial 
information, since it implies that $V$ is an integer, mod $4\pi^{2}/3.$ 
The renormalized volume in this case may serve as an analogue of 
Thurston's theory of the volume of closed hyperbolic 3-manifolds. 
(After completion of the paper, the referee informed the author that 
the equality in (0.5) for hyperbolic manifolds has also been proved 
by C. Epstein in Appendix A to [13]).

 A result analogous to (0.4) holds for AHE metrics on $M = M^{n}$ in 
any dimension $n \geq $ 2 and relates the Chern-Gauss-Bonnet integrand 
(Euler density) with $\chi (M)$ and the volume renormalization $V$ in 
even dimensions. In odd dimensions, it (re)-produces the formula for 
the conformal anomaly, c.f. [10]. This will be detailed elsewhere, and 
we restrict here to dimension 4.

\medskip

 It is interesting to compare, and combine, Theorem 0.1 with a result 
of Hitchin [11], where an analogous result is proved for the signature 
via the Atiyah-Patodi-Singer index theorem. Thus for any AH Einstein 
metric, (or any AH metric suitably asymptotic to an Einstein metric at 
infinity), one has
\begin{equation} \label{e0.6}
\frac{1}{12\pi^{2}}\int_{M}(|W^{+}|^{2} -  |W^{-}|^{2})dV = \tau (M) -  
\eta_{\gamma}, 
\end{equation}
where $\tau (M)$ is the signature of $M$ and $\eta_{\gamma}$ is the 
eta-invariant of the conformal infinity $(\partial M, \gamma ).$ In 
particular, (0.4) and (0.6) imply the following analogue of the 
Hitchin-Thorpe inequality for AHE metrics:
\begin{equation} \label{e0.7}
\chi (M) -  \frac{3}{4\pi^{2}}V \geq  \frac{3}{2}|\tau (M) -  
\eta_{\gamma}|, 
\end{equation}
with equality if and only if $(M, g)$ is self-dual Einstein.

 The volume term $V$ clearly depends, apriori, on the global geometry 
of the 'bulk' manifold $(M, g)$. However, the $\eta$-invariant of 
$(\partial M, \gamma )$ depends only on the intrinsic geometry of the 
conformal class $[\gamma ]$ on $\partial M.$ Thus, for a self-dual 
Einstein metric on $M$, it follows that $V$ is also an intrinsic 
invariant of $(\partial M, [\gamma ])$, given that the topology of $M$ 
is fixed. 

 The second purpose of this paper is to discuss to what extent this 
might be true for a general AH Einstein metric on $M$. First, recall 
that $\eta $ is a global conformal invariant of $[\gamma ],$ i.e. it is 
not computable from the local geometry of $[\gamma ].$ However, the 
variation of $\eta $ in the space of metrics is a local quantity; thus, 
if $h_{(0)}$ is an infinitesimal variation of $\gamma $ on $\partial 
M,$ then
\begin{equation} \label{e0.8}
d\eta (h_{(0)}) = -\frac{1}{24\pi^{2}}\int_{\partial M}< *dRic, 
h_{(0)}> dvol_{\gamma}, 
\end{equation}
c.f. [4, Thm. 6.9] and [1, Prop. 4.19]. Here $Ric$ is the Ricci 
curvature of $\gamma ,$ viewed as a 1-form with values in the tangent 
bundle, $d$ is the exterior derivative on $\Lambda^{1}$ induced by the 
metric, and $*$ is the Hodge star operator $*: \Lambda^{2} \rightarrow  
\Lambda^{1}.$ Recall that $dRic$ is the well-known Cotton-York tensor 
of conformal geometry, whose vanishing characterizes conformal flatness.

 Now let $g$ be an AH Einstein metric on $M$ and $h$ an infinitesimal 
AHE variation of $g$, with $h_{(0)}$ the induced variation of the 
boundary metric $\gamma .$ We then have the following expression for 
the variation or differential $dV$ of $V$.
\bbgin{theorem} \label{t 0.2.}
  Let $g$ be an AH Einstein metric, and $h$ an infinitesimal AHE 
deformation. Then the differential of the renormalized volume $V$ in 
the direction $h$ is given by
\begin{equation} \label{e0.9}
dV(h) = -\frac{1}{4}\int_{\partial M}<g_{(3)}, h_{(0)}> dvol_{\gamma}, 
\end{equation}
where $g_{(3)}$ is the $3^{rd}$ order term in the Taylor expansion of 
the metric $\bar g$ at $\partial M$, w.r.t. the special defining 
function determined by $\gamma$, c.f. \S 1.
\end{theorem}

  A formula similar to (0.9) holds in all dimensions $\geq 4$. Thus the 
variation of $V$ at $g$ is determined solely by the behavior of $\bar g$ 
at the boundary $\partial M$. Formally speaking, we may consider $g_{(3)}$ 
as the gradient of the volume function $V$, modulo the factor 
$-\frac{1}{4}$. The term $g_{(3)}$ is formally undetermined, in the sense 
that the Einstein equations do not determine any local 
expression for $g_{(3)}$ at $\partial M$, c.f. [6], [7]. 
This is in contrast to the situation for the terms $g_{(j)}, j \leq 2$, 
which are determined locally by the geometry of $\gamma = g_{(0)}$.

  In Proposition 2.6, we relate the formulas (0.8) and (0.9). Namely, let 
$g$ be an AH Einstein metric on $M$ with boundary metric $\gamma$. If 
$\gamma$ is not conformally flat and $dV \neq 0$, then
\begin{equation} \label{e0.10}
dV(h) = -\frac{1}{12}\int_{\partial M}< *dRic, h_{(0)}^{+}> 
dvol_{\gamma} +\frac{1}{12} \int_{\partial M}< *dRic, h_{(0)}^{-}> 
dvol_{\gamma}, 
\end{equation}
where $h_{(0)}^{\pm} = \pi^{\pm}(h_{(0)})$ and $\pi^{\pm}$ are linear 
projection operators on the space of symmetric bilinear forms 
$S^{2}(\partial M)$. As in (0.9), the projections $\pi^{\pm}$ depend, 
apriori, on the term $g_{(3)}$.

  It is a rather delicate open question whether the dependence of $dV$ 
in (0.9) or (0.10) on $g_{(3)}$ can be reduced to a dependence 
only on the boundary metric $(\partial M, \gamma)$, as is the case for 
$d\eta$. We point out at the end of \S 2 that if $dV$ does depend only on 
$(\partial M, \gamma)$, then it must also depend on the global topology 
of the bulk or filling manifold $M$, (again in contrast to $d\eta$, which 
is independent of $M$). This is illustrated by observing that the hyperbolic 
metric on $H^{4}(-1)/{\Bbb Z} \approx {\Bbb R}^{3} \times S^{1}$ and the 
Schwarzschild AdS metric on ${\Bbb R}^{2} \times S^{2}$ have 
$\partial M = S^{2} \times S^{1}$, and with suitable normalization, have 
the same conformally flat boundary metric $\gamma$ on $S^{2} \times S^{1}$. 
However, both $V$ and $dV$ are different for these metrics.

  On the positive side, we will show elsewhere that an AH Einstein metric 
on a given manifold $M$ is uniquely determined, up to diffeomorphism, by 
$dV$ and the boundary metric $\gamma$, at least when the induced map 
$\pi_{1}(\partial M) \rightarrow \pi_{1}(M)$ is surjective. In addition, 
the results of this paper will be applied elsewhere to study the 
Dirichlet problem for AHE metrics with prescribed conformal infinity, 
c.f. [6], [8].

 I would like to thank Jack Lee and Claude LeBrun for interesting 
conversations on these topics and the referee for comments on the 
manuscript.

\section{Chern-Gauss-Bonnet and $V$.}

\setcounter{equation}{0}

 This section is concerned with the proof of Theorem 0.1. Before 
starting the proof, we discuss some further background material on 
conformally compact metrics, c.f. also [7], [8].

 If $g$ is a complete conformally compact metric on $M$, with defining 
function $\rho ,$ let
\begin{equation} \label{e1.1}
r = log \rho^{-1}, \rho  = e^{- r}. 
\end{equation}
A simple computation shows that
\begin{equation} \label{e1.2}
|\nabla r|_{g}^{2} = |\nabla\rho|_{\bar g}^{2} \equiv  |\bar 
\nabla\rho|^{2}, 
\end{equation}
and that this quantity is independent of the choice of defining 
function $\rho $ at $\partial M$. Hence when the compactification $\bar 
g$ is $C^{1}$, $|\bar \nabla\rho|^{2}$ on $\partial M$ is an invariant 
of the conformal structure $(\partial M, [\gamma ]).$

 Now a computation for the change of curvature under conformal change 
in the metric shows
\begin{equation} \label{e1.3}
\bar K_{ij} = \rho^{-2}(K_{ij} +|\bar \nabla\rho|^{2}) -  
\rho^{-1}\{\bar D^{2}\rho (\bar e_{i}, \bar e_{i}) + \bar D^{2}\rho 
(\bar e_{j}, \bar e_{j})\}, 
\end{equation}
and, if $(M, g)$ is Einstein with $Ric = - 3g,$ then
\begin{equation} \label{e1.4}
\bar Ric = - 2\rho^{-1}\bar D^{2}\rho  + \{3\rho^{-2}(|\bar 
\nabla\rho|^{2} -  1) -  \rho^{-1}\bar \Delta \rho \}\bar g. 
\end{equation}
Here $\bar K_{ij},$ (resp. $K_{ij}),$ denotes the sectional curvature 
of $(M, \bar g),$ (resp. $(M, g)$), in the $(\bar e_{i}, \bar e_{j})$ 
direction, where $\{\bar e_{i}\}$ form an orthonormal basis w.r.t. $(M, 
\bar g).$ Hence if $\bar g$ is $C^{2}$ smooth up to $\partial M,$ then 
\begin{equation} \label{e1.5}
K_{ij} = -|\bar \nabla\rho|^{2} + O(\rho^{2}). 
\end{equation}

 Thus, the complete metric $g$ is asymptotically of variable strictly 
negative curvature; the curvature varies between two negative 
constants. The metric $g$ is called {\it  asymptotically hyperbolic}  
(AH) if the invariant $|\bar \nabla\rho|^{2}$ satisfies
\begin{equation} \label{e1.6}
|\bar \nabla\rho|^{2} = 1 \ \ {\rm on} \ \  \partial M. 
\end{equation}
Note that if $(M, g)$ is Einstein, then (1.4) implies that (1.6) must 
hold, so that any conformally compact Einstein metric is automatically 
AH.

 It is also elementary to see, (c.f. [8]) that if $(M, g)$ is AH, then 
there is a defining function $\rho $ such that in a collar neighborhood 
$U$ of $\partial M,$
\begin{equation} \label{e1.7}
|\nabla r| = |\bar \nabla\rho| \equiv  1, 
\end{equation}
in $U$. The metrics $g$ and $\bar g$ thus split in $U$ as
\begin{equation} \label{e1.8}
g = dr^{2} + g_{r} \ {\rm and} \  \bar g = d\rho^{2} + \bar g_{\rho}, 
\end{equation}
where $g_{r} = \rho^{-2}\cdot  g_{\rho}$ is a curve of metrics on the 
3-manifold $\partial M.$ Thus, w.r.t. the metric $g$ or $\bar g,$ the 
flow lines of $\nabla r$ or $\bar \nabla\rho $ are geodesics. The 
function $\rho $ gives the distance to $\partial M$ w.r.t. $\bar g,$ 
while the function $r =$ log $\rho^{-1}$ is a distance function w.r.t. 
$g$ from the boundary of some compact set in $M$.

 Defining functions satisfying (1.7) are called {\it  special}, or 
alternately {\it geodesic}, defining functions. Special defining 
functions are still not unique; as an example, $r$ may be the distance 
function from the boundary of any compact convex subset of 
$(M, g) = H^{4}(- 1).$ 

 If $(M, g)$ has a special $C^{k}$ conformal compactification $(M, \bar g),$ 
then one may expand $\bar g,$ i.e. $\bar g_{\rho}$ in (1.8), in a 
Taylor series in powers of $\rho ,$ as
\begin{equation} \label{e1.9}
\bar g_{\rho} = g_{(0)} + \rho g_{(1)} + \rho^{2}g_{(2)} + 
\rho^{3}g_{(3)} + ... + \rho^{k}g_{(k)} + o(\rho^{k}), 
\end{equation}
where the terms $g_{(i)}$ are bilinear forms on $T(\partial M),$ i.e. 
are annihilated when evaluated on $\bar \nabla \rho .$ The term 
$g_{(0)}$ is just the boundary metric $\gamma ,$ while $g_{(j)} = 
\frac{1}{j!}{\cal L}_{\bar \nabla \rho}^{(j)}\bar g_{\rho}|_{\rho =0},$ 
where ${\cal L} $ is the Lie derivative.

 Now if $g$ is Einstein, then results of Fefferman-Graham [6], c.f. 
also [7], imply that
\begin{equation} \label{e1.10}
g_{(1)} = 0, 
\end{equation}
so that $\partial M$ is totally geodesic in $\bar M$ w.r.t. $\bar g$ 
and further that
\begin{equation} \label{e1.11}
trg_{(3)} = 0, \ \ \delta g_{(3)} = 0, 
\end{equation}
where the trace and divergence are w.r.t. $\gamma .$ The term $g_{(2)}$ 
is intrinsically and locally determined by $\gamma  = g_{(0)},$ c.f. 
[7] or (2.18) below, but the Einstein equations do not imply any local 
intrinsic determination of $g_{(j)},$ for $j \geq $ 3, beyond (1.11).

 The expansion (1.9) gives the following expansion for $vol \bar S(\rho 
):$
\begin{equation} \label{e1.12}
vol \bar S(\rho ) = vol \bar S(0) + \rho^{2}v_{(2)} + O(\rho^{4}). 
\end{equation}
There is no $\rho^{3}$ term, by (1.11). Hence, $vol S(r)$ has the 
expansion
\begin{equation} \label{e1.13}
vol S(r) = v_{(0)}e^{3r}+ v_{(2)}e^{r} + O(e^{- r}), 
\end{equation}
and so, as in (0.3), $vol B(r)$ has the expansion
\begin{equation} \label{e1.14}
vol B(r) = \frac{1}{3}v_{(0)}e^{3r}+ v_{(2)}e^{r} + V + O(e^{-r}). 
\end{equation}

 We now begin with the proof of Theorem 0.1 itself. We assume that $g$ 
is an AHE metric on the 4-manifold $M$, and let $\rho $ be a special 
defining function for $(M, g)$. Since $g$ is Einstein, the curvature 
tensor $R$ is pure Weyl and scalar, i.e.
$$R = W -  \frac{1}{2}g \wedge g, $$
where $\wedge$ denotes the Kulkarni-Nomizu product, c.f. [2, Ch.1G]. 
When $W =$ 0, the curvature tensor $R = -\frac{1}{2}g \wedge g$ gives 
the curvature of hyperbolic space $H^{4}(- 1),$ with sectional 
curvature $- 1.$ Thus
\begin{equation} \label{e1.15}
\int_{M}(|R|^{2} -  6)dvol_{g} =\int_{M}|R+\frac{1}{2}g \wedge 
g|^{2}dvol_{g} = \int_{M}|W|^{2}dvol_{g} = \int_{M}|\bar 
W|^{2}dvol_{\bar g}, 
\end{equation}
where the second equality uses the conformal invariance of the $L^{2}$ 
norm of $W$ on 4-manifolds. Since, by assumption, $g$ has a $C^{2}$ 
conformal compactification, this integral is finite. The norm here is 
the usual $L^{2}$ norm of $R$ or $W$, as a symmetric map 
$\Lambda^{2}(TM) \rightarrow  \Lambda^{2}(TM),$ so that $|R|^{2} =$ 6 
on $H^{4}(- 1);$ this is $\frac{1}{4}|R|^{2},$ when $R$ is viewed as a 
(4,0) tensor.

 Let $D$ be a compact domain in $M$, with smooth boundary $\partial D 
\subset M$. Since $g$ is Einstein, the Chern-Gauss-Bonnet formula for 
manifolds with boundary states
\begin{equation} \label{e1.16}
\frac{1}{8\pi^{2}}\int_{D}|R|^{2} = \chi (D) -  
\frac{1}{2\pi^{2}}\int_{\partial D}\prod_{1}^{3}\lambda_{i} -  
\frac{1}{8\pi^{2}}\int_{\partial D}   \sum_{\sigma\in 
S_{3}}K_{\sigma_{1}\sigma_{2}}\lambda_{\sigma_{3}}, 
\end{equation}
c.f. [3]. Here $\lambda_{i}$ are the eigenvalues of the $2^{\rm nd}$ 
fundamental form $A$ and the indices $\sigma_{i}$ run over an 
orthonormal basis of the tangent spaces to $\partial D.$ The sign on 
$A$ is chosen so that $\lambda_{i} > $ 0 for convex domains; $K$ 
denotes sectional curvature.

 Let $\partial D = S(r)$ be the $r$-level set of the function $r$ in 
(1.1) and let $D = B(r)$ be the corresponding sublevel set. The $2^{\rm 
nd}$ fundamental form of $S(r)$ is then given by $A = D^{2}(r) = 
-D^{2}log \rho$. For $r$ sufficiently large, i.e. $\rho $ sufficiently 
small, $D$ is diffeomorphic to $M$. Hence (1.16) may be rewritten as
\begin{equation} \label{e1.17}
\frac{1}{8\pi^{2}}\int_{B(r)}|W|^{2} = \chi (M) -  \frac{3}{4\pi^{2}} 
\bigl ( vol B(r) + \frac{2}{3}\int_{S(r)}\prod_{1}^{3}\lambda_{i} + 
\frac{1}{6}\int_{S(r)} \sum_{\sigma \in S_{3}}K_{\sigma_{1} 
\sigma_{2}}\lambda_{\sigma_{3}} \bigr ) . 
\end{equation}

 All three terms in the parenthesis diverge to $\pm\infty $ as $r 
\rightarrow  \infty ,$ and so we need to understand their cancellation 
properties. From the expansion (1.9) and (1.14), using (1.3) and (1.18) 
below, one may prove purely formally that these terms must converge to 
$V$ as $r \rightarrow  \infty .$ However, it is worthwhile to calculate 
this explicitly to see just how the Einstein condition is being used. 
We begin by analysing the boundary integrals over $S(r)$. Following 
this, we analyse the bulk integral over $B(r)$.

 The eigenvalues $\lambda_{i}$ of $D^{2}r$ are related to the 
eigenvalues $\bar \lambda_{i}$ of $\bar D^{2}\rho $ by
\begin{equation} \label{e1.18}
\lambda_{i} = 1-\rho\bar \lambda_{i}, 
\end{equation}
c.f. [2, Ch. 1J] for instance for formulas on conformal changes of the 
metric. Hence
$$\prod\lambda_{i} = 1 -  \bar H\rho  + \bar \sigma\rho^{2} -  \bar 
\pi\rho^{3}, $$
where $\bar H =$ tr $\bar D^{2}\rho $ is the mean curvature of $S(r) = 
\bar S(\rho ), \bar \sigma = \prod_{i<j}\bar \lambda_{i}\bar 
\lambda_{j}$ and $\bar \pi = \bar \lambda_{1}\bar \lambda_{2}\bar 
\lambda_{3}.$ Here and in the following, the indices 1,2,3 refer to 
directions tangent to $S(r)$, while the index 4 refers to the normal 
direction.

 Next, for the boundary curvature term in (1.17), using (1.3) and 
(1.18) we have
$$\sum_{\sigma\in S_{3}}K_{\sigma_{1}\sigma_{2}}\lambda_{\sigma_{3}} = 
- 6 + 6\bar H\rho  + 2\bar \tau \rho^{2} + O(\bar \lambda^{2})\bar 
\rho^{2} + O(\bar K \bar \lambda)\rho^{3}, $$
where $\bar \tau = \bar K_{12}+\bar K_{13}+\bar K_{23},$ and $O(\bar 
\lambda^{2})$ and $O(\bar K \bar \lambda)$ denote terms quadratic in 
$\bar \lambda$ or products of $\bar K$ and $\bar \lambda$. For the last 
two terms in (1.17), we thus have
\begin{equation} \label{e1.19}
\int_{S(r)}\{(\frac{2}{3} -  1) + (-\frac{2}{3} + 1)\bar H\rho  + 
\frac{1}{3}\bar \tau\rho^{2} + O(\bar \lambda^{2})\rho^{2} + O(\bar 
\pi, \bar K \bar \lambda)\rho^{3}\}. 
\end{equation}
Now $vol_{g}S(r) = \rho^{-3}vol_{\bar g}S(\rho ) \sim  \rho^{- 3}.$ On 
the other hand, by (1.10), we have
\begin{equation} \label{e1.20}
\bar D^{2}\rho  = \bar A = \frac{1}{2}{\cal L}_{\bar \nabla \rho}\bar g 
= O(\rho ), 
\end{equation}
where $\bar A$ is the $2^{\rm nd}$ fundamental form of $\bar S(\rho )$ 
in $(M, \bar g).$ Hence, the last two terms in (1.19) are 
$O(\rho^{4}).$ This shows that (1.19) may be rewritten in the form
$$-\frac{1}{3}vol S(r) +  \frac{1}{3}\int_{S(r)}(\bar H\rho  + \bar 
\tau\rho^{2}) + O(\rho ). $$
We rewrite the second term as follows. From (1.4), one computes 
\begin{equation} \label{e1.21}
\bar Ric(4,4) = \frac{1}{6}\bar s = -\frac{\bar \Delta \rho}{\rho}, 
\end{equation}
where $\bar Ric$ denotes Ricci curvature w.r.t. $\bar g$. The first 
equality gives $5\bar Ric(4,4) = \sum_{i < 4}\bar Ric(i,i) = \bar Ric(4,4) + 
2\bar \tau,$ and so $2\bar Ric(4,4) = \bar \tau = - 2\frac{\bar 
H}{\rho}.$ Hence
$$\frac{1}{3}\int_{S(r)}(\bar H\rho  + \bar \tau\rho^{2}) = 
-\frac{1}{3}\int_{S(r)}\bar H\rho . $$
Finally, the integral curves of $\bar \nabla\rho $ are geodesics, and 
so the Ricatti equation
$$\frac{d\bar H}{d\rho} + |\bar A|^{2} + \bar Ric(\bar \nabla\rho , 
\bar \nabla\rho ) = 0, $$
holds. Since $|\bar A|^{2} = O(\rho^{2})$ by (1.20), and $\bar H/\rho = 
- \bar Ric(\bar \nabla\rho , \bar \nabla\rho )$, we obtain

$$\int_{S(r)}\bar H\rho  = \int_{S(r)}\bar H'\rho^{2} + O(\rho ), $$
where $\bar H'  = d\bar H/d\rho .$ In summary, we thus have the last 
two terms in (1.17) equal to
\begin{equation} \label{e1.22}
-\frac{1}{3}vol S(r) -  \frac{1}{3}\rho^{2}\int_{S(r)}\bar H'  + O(\rho 
). 
\end{equation}

 Now we claim that the two terms in (1.22) are exactly the first two 
terms in the $\rho$-expansion of $vol B(r)$.

\bbgin{lemma} \label{l 1.1.}
  As $r \rightarrow  \infty ,$ we have the expansion
\begin{equation} \label{e1.23}
vol B(r) = \frac{1}{3}vol S(r) + \frac{1}{3}\rho^{2}\int_{S(r)}\bar H'  
+ V + o(1). 
\end{equation}
\end{lemma}
{\bf Proof:}
 Let $\bar S(\rho )$ be the $\rho$-level set of $\rho $ in $(M, \bar 
g).$ Then for $\rho $ small,
$$vol \bar S(\rho ) = vol \bar S(0) + \rho\int_{\bar S(0)}\bar H + 
\frac{1}{2}\rho^{2}\int_{\bar S(0)}(\bar H' +\bar H^{2}) +  
\frac{1}{6}\rho^{3}\int_{\bar S(0)}(\bar H'' +3\bar H\bar H'  + \bar 
H^{3}) + O(\rho^{4}). $$
Recall that $vol S(r) = \rho^{-3}\cdot  vol \bar S(\rho )$ and $\bar H 
= \bar H''  =$ 0 at $\bar S(0) = \partial M$ by (1.10) and (1.11). Thus
\begin{equation} \label{e1.24}
vol S(r) = \rho^{-3}vol \bar S(0) + \frac{1}{2}\rho^{- 1}\int_{\bar 
S(0)}\bar H'   + O(\rho ). 
\end{equation}
Integrating this from 0 to $r$ gives
$$vol B(r) = \int_{0}^{r}vol S(r)dr = \int_{\rho}^{1}\rho^{-1}vol \bar 
S(\rho)d\rho  = vol \bar S(0)\int_{\rho}^{1}\rho^{-4}d\rho  + 
(\frac{1}{2}\int_{\bar S(0)}\bar H' )\int_{\rho}^{1}\rho^{-2}d\rho  + 
O(1), $$
which implies
$$vol B(r) = \frac{1}{3}\rho^{-3}vol \bar S(0) + 
\frac{1}{2}\rho^{-1}\int_{\bar S(0)}\bar H'  + V + o(1). $$
Substituting in (1.24) shows that
\begin{equation} \label{e1.25}
vol B(r) = \frac{1}{3}vol S(r) + \frac{1}{3}\rho^{-1}\int_{\bar 
S(0)}\bar H'  + V + o(1). 
\end{equation}
Finally, we have
$$\rho^{-1}\int_{\bar S(0)}\bar H'  = \rho^{-1}\int_{\bar S(\rho )}\bar 
H' (\rho ) + \rho^{-1}(\int_{\bar S(0)}\bar H'  -\int_{\bar S(\rho 
)}\bar H' (\rho )).$$ 
But $\bar H' (\rho ) = \bar H' (0) + \rho\bar H'' (0) + o(\rho^{2}) = 
\bar H' (\rho ) + o(\rho^{2}).$ Hence
$$\rho^{-1}\int_{\bar S(0)}\bar H'  = \rho^{-1}\int_{\bar S(\rho )}\bar 
H' (\rho ) + o(1),$$ 
and the result follows.

{\endproof}

 Combining (1.17), (1.22) and (1.23) and letting $r \rightarrow  \infty 
$ then completes the proof of Theorem 0.1.
\bbgin{remark} \label{r 1.2.}
  {\rm The same proof as above evaluates the right side of (1.17) 
whenever $(M, g)$ is any AH metric which is Einstein to $3^{\rm rd}$ 
order, i.e. for which the expansion (1.9) agrees with the expansion of 
an Einstein metric to order 3. Hence, for such metrics, we obtain}
\begin{equation} \label{e1.26}
\frac{1}{8\pi^{2}}\int_{M}(|W|^{2} -  \frac{1}{2}|z|^{2} + 
\frac{1}{24}s^{2}- 6)dV = \chi (M) -  \frac{3}{4\pi^{2}}V. 
\end{equation}

 {\rm It follows for instance that an Einstein metric minimizes $V$ in 
its conformal class, among AH metrics. Note that $V$ itself is, of 
course, not a conformal invariant among such AH metrics.}
\end{remark}

\section{Boundary Determination of $dV$.}

\setcounter{equation}{0}

 This section is concerned with the question of to what extent the 
renormalized volume $V$, or the $L^{2}$ norm of the Weyl curvature, is 
determined by the conformal infinity $\gamma $ of an AHE metric. To do 
this, we study the variation $dV$ of $V$ in the space of AHE metrics on 
$M$.

 Thus, let $g$ be an AHE metric on $M$ and let $h$ be an infinitesimal 
variation of $g$, so that the curve of metrics $g_{t} = g+th$ is AHE, 
to first order in $t$. From Theorem 0.1, we have
\begin{equation} \label{e2.1}
dV(h) = \frac{dV_{t}}{dt}|_{t=0} = 
-\frac{1}{6}\frac{d}{dt}(\int_{M}|W|^{2}dvol_{t})_{t=0} \equiv 
-\frac{1}{6}d{\cal W} (h). 
\end{equation}
To analyse $dV$ recall that, by definition, Einstein metrics are 
critical points of the scale-invariant Einstein-Hilbert action
\begin{equation} \label{e2.2}
{\cal S}  = vol^{-1/2}\int sdvol, 
\end{equation}
in dimension 4. Hence the variation $d{\cal S} $ of ${\cal S} $ is 
determined by the behavior of the variation of the metric at the 
boundary. We first make this precise in the Lemma below, and then 
relate it to the variation of $V$. The following result has recently 
also been proved in [15]; the proof below however is simpler and more 
transparent.

\bbgin{lemma} \label{l 2.1.}
  Let $g$ be an Einstein metric on a smooth compact domain $D$ in 
$M^{4}$, with scalar curvature s, and let $h$ be an infinitesimal 
deformation of g, so that $g_{t} = g+th$ is Einstein, with scalar 
curvature s, to first order in t. Then
\begin{equation} \label{e2.3}
(vol D)'  = \frac{d}{dt}vol_{g_{t}}D|_{t=0} = 
-\frac{2}{s}\int_{\partial D}(2H'  + < A, h> ), 
\end{equation}
where $A$ is the $2^{nd}$ fundamental form of $\partial D, H = tr A$, 
and $H'  = \frac{dH}{dt}|_{t=0.}$
\end{lemma}
{\bf Proof:}
 Take the derivative of (2.2) w.r.t. $t$ and use the fact that $s$ is 
constant. A brief computation shows that at $t = 0$,
$$s\cdot  (vol^{1/2}D)'  = vol^{-1/2}\int_{D}(L(h) + \frac{s}{4}<g, h> 
)dvol, $$
where $L(h) = s' (h)$ is the linearization of the scalar curvature, 
given by
$$L(h) = -\Delta trh + \delta\delta h \ -  < Ric, h> , $$
c.f. [2, 1.174]. Since $z =$ 0, this gives
$$\frac{1}{2}s(vol D)'  = \int_{D}(-\Delta trh + \delta\delta h)dvol, $$
and hence by the divergence theorem
\begin{equation} \label{e2.4}
\frac{1}{2}s(vol D)'  = -\int_{\partial D}< dtrh, N>  -  \int_{\partial 
D}\delta h(N), 
\end{equation}
where $N$ is the unit outward normal. 

 Now choose local normal exponential (Fermi) coordinates for a 
neighborhood of $\partial D.$ Thus, $N$ is the field tangent to 
geodesics, and normal to equidistant hypersurfaces $S(r)$, with $S(0) = 
\partial D.$ Let $\{e_{i}\}$ be a local orthonormal basis for 
$T(S(r))$, so that $\{e_{i}, N\}$ are a local orthonormal basis for 
$T(D)$ near $\partial D.$ We then have
\begin{equation} \label{e2.5}
\delta h(N)  = - Nh(N,N) \ -  <\nabla_{e_{i}} h(e_{i}), N> , 
\end{equation}
and
\begin{equation} \label{e2.6}
<dtrh, N> \ = Nh(N,N) + N(< g^{T}, h^{T}>).  
\end{equation}
where $T$ denotes tangential part. When combined, the first terms in 
(2.5) and (2.6) cancel. For the second term in (2.5), we have 
$<\nabla_{e_{i}}h(e_{i}), N>  = div_{S}(h(N)) -  < h, A> ,$ where 
$div_{S}$ is the divergence on the hypersurfaces $S$. This integrates 
to 0 on $S(0) = \partial D.$ Hence, (2.4) becomes
\begin{equation} \label{e2.7}
\frac{1}{2}s(vol D)'  = -\int_{\partial D}N< g^{T}, h^{T}>  -  
\int_{\partial D}< A, h> . 
\end{equation}
To evaluate the first term, for each metric $g_{t}$ we have the 
hypersurfaces $S_{t}(r)$ constructed above, with the induced metric 
$g_{ij}(t, r)$. In a fixed local coordinate system, the volume form 
$dV_{S}(t, r)$ of $S_{t}(r)$ is given by
$$dV_{S}(t,r) = (det g_{ij}(t,r))^{1/2}dx_{1}\wedge dx_{2}\wedge 
dx_{3}. $$
Then $\frac{1}{2}< g^{T}, h^{T}>dV_{S} = \frac{\partial}{\partial 
t}[det g_{ij}(t,r)^{1/2}]dx_{1}\wedge dx_{2}\wedge dx_{3},$ and
$$\tfrac{1}{2}N< g^{T}, h^{T}>dV_{S} = \frac{\partial}{\partial 
r}\frac{\partial}{\partial t}[det g_{ij}(t,r)^{1/2}]dx_{1}\wedge 
dx_{2}\wedge dx_{3}. $$
The coefficients $g_{ij}$ are smooth functions of the parameters $r$ 
and $t$, and so
$$\frac{\partial}{\partial r}(\frac{\partial}{\partial 
t}(det(g_{ij})^{1/2})) = \frac{\partial}{\partial 
t}(\frac{\partial}{\partial r}(det(g_{ij})^{1/2})) = H'  = 
\frac{dH}{dt}. $$
It follows that (2.7) becomes
$$\frac{1}{2}s(vol D)'  = -\int_{\partial D}(2H'  + < A, h> )dvol, $$
which gives (2.3).

{\endproof}

 This result, with the same proof, holds in all dimensions, with the 
coefficient $\frac{1}{2}$ replaced by $2/n, n =$ dim $M$.

 We now apply Lemma 2.1 to the domains $B(r)$ in an AH Einstein 
manifold $(M, g)$, and let $r \rightarrow  \infty .$ This leads to the 
proof of Theorem 0.2, which we restate as:
\bbgin{theorem} \label{t 2.2.}
  Let $h = dg_{t}/dt$ be an infinitesimal AHE deformation of an AHE 
metric $g$ on $M$ and let $h_{(0)} = d\gamma_{t}/dt$ be the induced 
variation of $\gamma $ on $\partial M,$ where $\bar g$ is determined 
by a special defining function $\rho ,$ as in (1.7). Then
\begin{equation} \label{e2.8}
dV(h) = -\frac{1}{4}\int_{\partial M}< g_{(3)}, h_{(0)}> , 
\end{equation}
for $g_{(3)}$ as in (1.9). The inner product and volume form in (2.8) 
are w.r.t. $\gamma .$
\end{theorem}
{\bf Proof:}
 By Lemma 2.1, we have with $r = log \rho^{-1} >> $ 1,
\begin{equation} \label{e2.9}
(vol B(r))'  = \frac{1}{6}\int_{S(r)}(2H'  + < A, h> )dvol = 
\frac{1}{6}\rho^{-3}\int_{\bar S(\rho )}(2H'  + < A, h> )d\bar vol. 
\end{equation}
As in the proof of Theorem 0.1, we analyse the terms on the right from 
the expansion of $\bar g_{t},$ given by
\begin{equation} \label{e2.10}
\bar g_{t} = d\rho_{t}\otimes d\rho_{t} + (g_{(0),t} + 
\rho_{t}^{2}g_{(2),t} + \rho_{t}^{3}g_{(3),t}). 
\end{equation}
Taking the derivative of (2.10) w.r.t. $t$ gives
\begin{equation} \label{e2.11}
\bar h = 2d\rho'\otimes d\rho  + (h_{(0)} + \rho^{2}h_{(2)} + 
2\rho\rho' g_{(2)} + O(\rho^{3})). 
\end{equation}
Here $\rho'  = d\rho_{t}/dt,$ and we have used the fact that $\rho'  = 
O(\rho ),$ since $\rho' (\partial M) =$ 0. In fact, since $\rho_{t}$ 
are special defining functions w.r.t. $g_{t},$ a simple computation 
gives 
\begin{equation} \label{e2.12}
\rho'  = \phi_{(1)}\rho  + \phi_{(3)}\rho^{3} + o(\rho^{3}), 
\end{equation}
c.f. also [7, Lemma 2.2]. Note also that $\bar h$ has no tangential 
terms of order $\rho .$

 Next, from (1.18), we have $H = 3-\bar H\rho ,$ and, by (1.10)-(1.11), 
$\bar H = h_{1}\rho  + h_{3}\rho^{3},$ so that
\begin{equation} \label{e2.13}
H = 3- h_{1}\rho^{2} -  h_{3}\rho^{4}. 
\end{equation}
Similarly, $A = g^{T} -  \rho\bar A,$ and $\bar A = A_{1}\rho  + 
A_{2}\rho^{2},$ so that
\begin{equation} \label{e2.14}
A = g^{T} -  A_{1}\rho^{2} -  A_{2}\rho^{3} + O(\rho^{4}). 
\end{equation}
Further, by (2.13), $H'  = h_{1}'\rho^{2} + 2\rho\rho' h_{1} + 
O(\rho^{4})$ which with (2.12) gives
\begin{equation} \label{e2.15}
H'  = \xi_{1}\rho^{2} + O(\rho^{4}). 
\end{equation}

 Now we substitute these computations in (2.9). The estimate (2.15) 
shows that the first term in (2.9) contains only an $O(\rho^{-1})$ 
term, and hence gives no contribution to $V' ,$ where $V$ is the 
renormalized volume. Hence we may ignore this term. For the next term, 
we have
\begin{equation} \label{e2.16}
< A, h> \ = \ < h, g^{T}>  -  <A_{1}, h>\rho^{2} -  < A_{2}, h>\rho^{3} 
+ O(\rho^{4}). 
\end{equation}
Now for any (1,1) tensors $A, B$, $< A, B>_{g} \ = \ <A, B>_{\bar g},$ 
and so the $g$-inner products in (2.16) may be replaced by $\bar 
g$-inner products. We have $\frac{1}{2}<h, g^{T}>_{\bar g}dvol_{\bar 
S(\rho )} = \frac{d}{dt}(dvol_{\bar g_{t}}(\bar S_{t}(\rho ))),$ which 
vanishes at order $O(\rho^{3})$ by (1.11). Similarly, $<A_{1}, h>_{\bar 
g}\rho^{2}$ has no terms of order $O(\rho^{3})$ by (2.11). Hence the 
only term in (2.16) of order $O(\rho^{3})$ is
$$-< A_{2}, h>\rho^{3} = -<A_{2}, h>_{\bar g}\rho^{3} = -<A_{2}, 
h_{(0)}>_{\bar g}\rho^{3} + o(\rho^{3}).$$ 
Taking the limit $\rho  \rightarrow $ 0 then implies that
$$V'  = -\frac{1}{6}\int_{\partial M}< A_{2}, h_{(0)}> .$$ 
To complete the proof, we have ${\cal L}_{\bar \nabla\rho}\bar g = 
2\bar A$, while $A_{2} = \frac{1}{2}{\cal L}_{\bar \nabla\rho}^{2}\bar 
A = \frac{1}{4}{\cal L}_{\bar \nabla\rho}^{3}\bar g = 
\frac{3}{2}g_{(3)},$ which then gives (2.8).

{\endproof}

 The formula (2.8) shows that although apriori the renormalized volume 
$V$ depends on the global geometry of the bulk manifold $(M, g)$, its 
variation $dV$ depends only on the $(3^{\rm rd}$ order) behavior of the 
compactification $\bar g$ at $\partial M.$ We note that one may prove, 
via Lemma 2.1 again, that there is a similar formula in higher 
dimensions.
\begin{remark} \label{r 2.3.}
  {\rm Using the fact that $g_{(2)} = {\cal L}_{\bar \nabla\rho}\bar A$ 
together with (1.4), a brief computation shows that}
\begin{equation} \label{e2.17}
g_{(2)} = -\frac{1}{2}(\bar Ric -  \frac{\bar s}{6}g_{(0)}), 
\end{equation}
{\rm where the curvatures are w.r.t. the metric $\bar g$ on $M$. In 
fact one may compute that}
\begin{equation} \label{e2.18}
-\frac{1}{2}(\bar Ric -  \frac{\bar s}{6}g_{(0)}) = - (Ric_{\gamma} -  
\frac{s_{\gamma}}{4}\gamma ), 
\end{equation}
{\rm at $\partial M,$ where the curvatures on the right of (2.18) are 
{\it intrinsic}  w.r.t. the boundary metric $g_{(0)} = \gamma ,$ c.f. 
also [7, (2.10)] for example. We note that (1.21) and (1.11) imply that 
$\bar \nabla\rho (\bar s) =$ 0 at $\partial M$ and further computation 
shows that $(\bar \nabla_{X}\bar Ric)(\bar \nabla\rho ) =$ 0 at 
$\partial M$, for $X$ tangent to $\partial M$. Hence, (2.17) and the 
relation $g_{(3)} = \frac{1}{3}{\cal L}_{\bar \nabla\rho}g_{(2)}$ imply}
\begin{equation} \label{e2.19}
g_{(3)} = \frac{1}{6}\nabla_{N}\bar Ric = \frac{1}{6}d\bar Ric(N), 
\end{equation}
{\rm where $N = -\bar \nabla\rho $ and $d = d^{\bar \nabla}$ is the 
exterior derivative w.r.t. the ambient metric $\bar g,$ c.f. [2, 4.69]}.
\end{remark}

\begin{remark} \label{r 2.4.}
  {\rm We verify briefly that (2.8) also gives, up to a constant, the 
variation $d{\cal W} ,$ where ${\cal W}  = \int|W|^{2};$ of course, 
this must be the case by (2.1). 

 The gradient $\nabla{\cal W} ,$ i.e the Euler-Lagrange operator for 
${\cal W} ,$ is given by the Bach tensor $\nabla {\cal W}^{2} = 
\frac{1}{2}(\delta\delta W +  W\circ Ric),$ c.f. [2, 4.77]; (the factor 
of $\frac{1}{2}$ comes from the definition of $|W|^{2}$ as in \S 1). 
Einstein metrics are also critical points of ${\cal W} ,$ so that 
$\nabla{\cal W}^{2} =$ 0. Hence, as in the proof of Lemma 2.1, on any 
compact domain $D \subset  M$ and at $t =$ 0, we have
$$\frac{d}{dt}\int_{D}|W|^{2}dvol_{t} = \int_{D}\frac{d}{dt}(|W|^{2} 
dvol_{t}) = \int_{D}<\nabla{\cal W}^{2} , h>  + \int_{\partial D}< 
B_{W}, h> , $$
i.e.}
\begin{equation} \label{e2.20}
d{\cal W} (h) = \int_{\partial D}< B_{W}, h> , 
\end{equation}
{\rm where $B_{W}$ is a boundary term. To determine $B_{W},$ integrate 
by parts as follows:
$$\int_{D}<\delta\delta W, h>  = -\int_{\partial D}<\delta W(N), h>  + 
\int_{D}<\delta W, Dh> , $$
and
$$\int_{D}<\delta W, Dh> = -  \int_{\partial D}<W(N), Dh>  + \int_{D}< 
W, DDh> , $$
where $N$ is the unit outward normal. Hence 
\begin{equation} \label{e2.21}
\int_{\partial D}< B_{W}, h>  = -\frac{1}{2}\int_{\partial D}<\delta 
W(N), h>  -\frac{1}{2}\int_{\partial D}< Dh, W(N)> . 
\end{equation}
Since ${\cal W} $ is conformally invariant, we may compute (2.21) 
w.r.t. the compactification $\bar g$ and let $D =$ $M$. A computation 
as in that giving (2.19) shows that the second integral in (2.21) 
vanishes, and so
\begin{equation} \label{e2.22}
\int_{\partial M}< B_{W}, h_{(0)}>  = -\frac{1}{2}\int_{\partial 
M}<\delta W(N), h_{(0)}>  = \frac{1}{4}\int_{\partial M}< d(\bar Ric -  
\frac{\bar s}{6}\bar g)(N), h_{(0)}> , 
\end{equation}
where the second equality uses the Bianchi identity, c.f. [2, 16.3]. 
Via (2.8) and (2.19), this confirms (2.1).}

\end{remark}

\medskip

 Combining (2.8) and (2.19), we have
\begin{equation} \label{e2.23}
dV(h) = -\frac{1}{24}\int_{\partial M}< d\bar Ric(N), h_{(0)}> . 
\end{equation}
This formula resembles, at least formally, the formula for the 
variation of $\eta $ in (0.8), i.e.
\begin{equation} \label{e2.24}
d\eta (h) = -\frac{1}{24\pi^{2}}\int_{\partial M} < *dRic, h_{(0)}> . 
\end{equation}

Of course by (0.7) these formulas must agree if $(M, g)$ is self-dual 
Einstein, and $h$ is an infinitesimal variation of such metrics. 
However, in (2.24), the Ricci curvature $Ric$ and exterior derivative 
$d$ are {\it intrinsic}, i.e. computed on the 3-manifold $\partial M$ 
w.r.t. the boundary metric $\gamma .$ On the other hand, in (2.23), the 
Ricci curvature $\bar Ric$ and $d$ are {\it  extrinsic}, computed w.r.t 
the ambient metric $\bar g$ at $\partial M.$

 As in (2.18), the term $g_{(2)}$ in the expansion (1.9) of $\bar g$ at 
$\partial M$ is local and intrinsic; recall that $g_{(1)} =$ 0. 
However, as pointed out in [6], the Einstein equations only imply the 
relations (1.11) on the third term $g_{(3)}$ in the expansion; the 
remaining parts of $g_{(3)}$ are formally undetermined. The term 
$g_{(3)}$ hence (may) depend on the bulk metric $\bar g$ at $\partial 
M,$ and not only on $\gamma .$ 

 In general, this issue is related to the unique solvability of the 
Dirichlet problem for AH Einstein metrics with prescribed conformal 
infinity. Namely, if, given a boundary metric $\gamma ,$ there is a 
unique AHE metric $g$ with conformal infinity $\gamma ,$ then the 
expansion terms $g_{(k)}$ are all necessarily determined, in some 
manner, by the intrinsic geometry of $\gamma  = g_{(0)}.$ On the other 
hand, if this is not the case, then some $g_{(k)}$ may not be 
determined from $\gamma .$

 Before proceeding further, we make several remarks.

\bbgin{remark} \label{r 2.5.}
  {\rm (i). We observe that
\begin{equation} \label{e2.25}
dV \neq  c \cdot d\eta  
\end{equation}
for any constant $c$, in general. This essentially follows from Theorem 0.1 
and a remark of Hitchin in [11]. Thus let $\gamma_{o}$ be the canonical 
round metric on $S^{3}$ and let $\gamma $ be any metric sufficiently close 
to $\gamma_{o}$ which is invariant under an orientation reversing 
reflection of $S^{3}.$ Since $\eta $ changes sign under orientation 
reversal, $\eta (\gamma ) =$ 0. The Graham-Lee theorem [8] shows that 
any such $\gamma $ may be filled in with an AH Einstein metric $g$, 
with $\gamma $ as conformal infinity. If we now take a curve of such 
metrics $g_{t}$ with boundary values $\gamma_{t},$ then $d\eta 
(\gamma_{t})/dt = d\eta (h) =$ 0, for all $t$. However, such curves 
$g_{t}$ will satisfy, for $t > $ 0, $V(g_{t}) \neq  V(g_{o}) = 
\frac{4\pi^{2}}{3}$ unless the curve is a constant curve, c.f. Theorem 
0.1. Hence, for some $t \neq $ 0 small, $dV(h) \neq  0$, which gives 
(2.25).

 We recall that $dRic$ is the {\it  only} local conformal invariant 
constructed from the metric in dimension 3 and hence $dV$ cannot be a 
locally defined intrinsic invariant of $\gamma $ in general.

 (ii). In [6], Fefferman and Graham consider formally the class of AH 
Einstein metrics for which $\bar g$ has an {\it  even}  expansion 
(1.9), i.e. $g_{(odd)} =$ 0. Of course, for such metrics Theorem 2.2 
gives 
$$dV = -\frac{1}{6}d{\cal W}  = 0. $$
Thus, such points are critical points of $V$ or ${\cal W} $ and so one 
would expect that there are very few such metrics, in the space of all 
AHE metrics.

 Any hyperbolic metric is even in this sense, since the metrics $\bar 
g_{\rho}$ in (1.8) are given by $\bar g_{\rho} = (1-\rho^{2})^{2}\cdot  
g_{(0)}.$ But Theorem 0.1 implies that any hyperbolic metric gives the 
maximal value for $V$ on $M$ and so of course this metric must be a 
critical point of $V$.}
\end{remark}

 Next we derive the formula (0.10). For a given AH Einstein metric 
$(M, g)$, let $W^{+}$ and $W^{-}$ be the self-dual and anti-self-dual 
parts of the Weyl curvature, so that $W = W^{+} + W^{-};$ (recall that 
$M$ is an oriented 4-manifold). As in Remark 2.4, we have $\nabla{\cal 
W}  = \frac{1}{2}(\delta\delta W +  W \circ Ric) =$ 0, and hence
$$\nabla{\cal W}^{+} = \frac{1}{2}(\delta\delta W^{+} + W^{+} \circ 
Ric) = 0, \ {\rm and} \ \nabla{\cal W}^{-} = \frac{1}{2}(\delta\delta 
W^{-} + W^{-} \circ Ric) = 0, $$
where ${\cal W}^{\pm}$ are the functionals $\int|W^{\pm}|^{2}.$ 

 The same reasoning as in (2.20) then shows that, for an infinitesimal 
AHE deformation $h$,
\begin{equation} \label{e2.26}
d{\cal W}^{\pm}(h) = \frac{d}{dt} \bigl( \int_{M}|W^{\pm}|^{2} \bigr)_{t=0} 
= \int_{\partial M}< B^{\pm}, h_{(0)}> , 
\end{equation}
where $B^{\pm}$ are boundary terms. These terms may be computed in the 
same way as $B_{W}$ in (2.22). Thus, informally, we may think of 
$\nabla{\cal W}^{\pm} = B^{\pm}.$ 

 Now consider the moduli space ${\cal M} $ of AHE metrics on $M$ which 
admit a $C^{3}$ conformal compactification $\bar g.$ The boundary 
values $\gamma $ of such metrics give a space ${\cal B} $ of metrics on 
$\partial M.$ The structure of ${\cal M} $ and ${\cal B} $ is not of 
concern here. Instead, we consider only the (formal) tangent spaces 
$T_{g}{\cal M} ,$ i.e. the vector space of solutions to the linearized 
AHE equations at a given $\gamma\in{\cal M}$, c.f. [2, Ch. 12] for 
background on linearized Einstein equations. Any infinitesimal AHE 
deformation $h\in T_{g}{\cal M} $ induces an infinitesimal deformation 
$h_{(0)}\in T_{\gamma}{\cal B} .$

 Let $T^{\pm}$ be the subspace of tangent vectors which leave ${\cal 
W}^{\mp}$ unchanged, to first order; thus
\begin{equation} \label{e2.27}
T^{+} = ker d{\cal W}^{-}, \ \ {\rm and} \ \ T^{-} = ker d{\cal W}^{+}. 
\end{equation}
Via (2.26), $T^{\pm}$ may also be viewed as subspaces of 
$T_{\gamma}{\cal B} $ orthogonal to $B^{\mp}.$ They are codimension 
1 hyperplanes of $T_{\gamma}{\cal B} $, except when $d{\cal W}^{+} = 0$ 
or $d{\cal W}^{-} = 0$. Of course $T^{+}\cap T^{-}$ consists of the 
variations which change neither ${\cal W}^{+}$ or ${\cal W}^{-}$ 
to first order. Observe that if $d{\cal W}^{+} = \lambda d{\cal W}^{-}$ 
for some $\lambda \neq 0$, then since a change in the orientation 
interchanges $d{\cal W}^{+}$ and $d{\cal W}^{-}$, we must have 
$\lambda^{2} = 1$, and so either $d{\cal W} = 0$, (when $\lambda = 1$), 
or $d\eta = 0$, (when $\lambda = 1$, by (0.6) and (0.8)).

  Now suppose that $d{\cal W}^{+}$ and $d{\cal W}^{-}$ are linearly 
independent. Then any tangent vector $h_{(0)}\in T_{\gamma}{\cal B}$ 
may be decomposed uniquely as
\begin{equation} \label{e2.28}
h_{(0)} = h_{(0)}^{+} + h_{(0)}^{0} + h_{(0)}^{-}, 
\end{equation}
where $h_{(0)}^{\pm} \in T^{\pm}$ and $h_{(0)}^{0}\in T^{+}\cap T^{-}$. 
If $d{\cal W}^{-} = 0$, we set $h^{-} = 0$ in (2.28) and similarly 
$h^{+} = 0$ if $d{\cal W}^{+} = 0$. The decomposition (2.28) defines the 
projection operators $\pi^{\pm}$ in (0.10). A similar decomposition holds 
for $h \in T_{g}{\cal M}$.

 The following result relates $d{\cal W}$ and $d\eta$, and gives (0.10) 
via (2.1).

\bbgin{proposition} \label{p 2.6.}
  Let $g$ be an AH Einstein metric on $M$, with boundary metric $\gamma$. 
Suppose $d{\cal W} \neq 0$ and $d\eta \neq 0$, i.e. $g$ is not a critical 
point of ${\cal W}$ or $\eta$. Then, in the notation above, we have
\begin{equation} \label{e2.29}
d{\cal W} (h) = \frac{1}{2}\int_{\partial M}< *dRic, h_{(0)}^{+}>  -  
\frac{1}{2}\int_{\partial M}< *dRic, h_{(0)}^{-}>, 
\end{equation}
where $*dRic$ is intrinsically defined w.r.t. $\gamma$.
\end{proposition}
{\bf Proof:}
 Since ${\cal W}  = {\cal W}^{+} + {\cal W}^{-},$ we have
$$d{\cal W} (h) = d{\cal W}^{+}(h) + d{\cal W}^{-}(h). $$
By the construction in (2.27)-(2.28), this gives
$$d{\cal W} (h) = d{\cal W}^{+}(h^{+}) + d{\cal W}^{-}(h^{-}), $$
where
$$d{\cal W}^{\pm}(h^{\pm}) = \int_{\partial M}< B^{\pm}, h_{(0)}^{\pm}> 
. $$
However since $h^{+}$ leaves ${\cal W}^{-}$ unchanged to first order, 
it follows from (0.6) and (0.8) that
\begin{equation} \label{e2.30}
d{\cal W}^{+}(h^{+}) = d({\cal W}^{+}(h^{+}) -  {\cal W}^{-}(h^{+})) = 
\frac{1}{2}\int_{\partial M}< *dRic, h_{(0)}^{+}> . 
\end{equation}
Thus, when viewed as a linear functional restricted to $T^{+}\subset  
T_{\gamma}{\cal B} , B^{+} = *dRic$ is intrinsically and locally 
determined by $\gamma .$ For the same reasons, we also obtain
\begin{equation} \label{e2.31}
d{\cal W}^{-}(h^{-}) = -\frac{1}{2}\int_{\partial M}< *dRic, 
h_{(0)}^{-}> . 
\end{equation}
Combining (2.30) and (2.31) gives the result.

{\endproof}

  While $*dRic$ is intrinsically determined by the boundary metric 
$\gamma$, it is not clear to what extent the subspaces $T^{\pm}$ are 
determined by $\gamma$.  In this regard, we discuss the following 
examples, which show that $d{\cal W}$ or $dV$ cannot be {\it solely} 
determined by the boundary metric, at least when the topological 
type of $M$ is allowed to vary; these examples are also discussed in 
[9] and [16].
   Thus, consider first hyperbolic 4-space ${\Bbb H}^{4}(-1)$. For any 
geodesic $\sigma \subset {\Bbb H}^{4}(-1)$, translation by a fixed 
length $L$ along $\sigma$ extends to an isometry of ${\Bbb H}^{4}(-1)$. 
Let ${\Bbb H}^{4}(-1)/{\Bbb Z} \approx {\Bbb R}^{3} \times S^{1}$ be the 
quotient, with ${\Bbb Z}$ the group generated by the translation. The 
metric $g_{-1}$ on ${\Bbb H}^{4}(-1)/{\Bbb Z}$ may be written as
$$g_{-1} = dr^{2} + sinh^{2}r g_{S^{2}(1)} + cosh^{2}d\theta^{2},$$
where $\theta$ parametrizes a circle of length $L$. We have 
$\partial M = S^{2} \times S^{1}$ and the conformal infinity $[\gamma]$ 
is the conformal class of the product metric $S^{2}(1) \times S^{1}(L)$. 
For instance by (0.4) and Remark 2.5(ii), we have $V = 0$ and 
$dV = g_{(3)} = 0$ for $g_{-1}$, and of course $\eta = d\eta = 0$.

  On the other hand, the AdS Schwarzschild metric, c.f. [9], [16, \S 3.2], 
or [2,9.118(d)] is an AH Einstein metric on ${\Bbb R}^{2} \times S^{2}$ 
given by
$$g_{AS} = (1+r^{2} - \tfrac{2m}{r})^{-1}dr^{2} + r^{2}g_{S^{2}(1)} + 
(1+r^{2} - \tfrac{2m}{r})d\theta^{2}.$$
Here $m > 0$ is the mass parameter, $r \geq r_{+}$, where $r_{+}$ is the 
largest root of the equation $1+r^{2} - \frac{2m}{r} = 0$, and $\theta$ 
parametrizes a circle of length $L = L(m) = \frac{4\pi r_{+}}{1+3r_{+}^{2}}$. 
This metric has the same conformal infinity $S^{2}(1) \times S^{1}(L)$ as 
before and so $\eta = d\eta = 0$. However, by [9, (2.9)],
$$V(g_{AS}) = \frac{\pi r_{+}^{2}(1-r_{+}^{2})}{1+3r_{+}^{2}},$$ 
and a straightforward computation using (0.9) gives 
$$dV_{g_{AS}} = \frac{m}{2}d\theta^{2}.$$

\bibliographystyle{plain}

\bigskip

\begin{center}
November, 2000
\end{center}

\medskip

\noindent
\address{Department of Mathematics \\
S.U.N.Y. at Stony Brook \\
Stony Brook, N.Y. 11794-3651}\\
\email{anderson@@math.sunysb.edu}

\end{document}